\documentclass[runningheads]{llncs}
\usepackage{amsmath, amsfonts}

\def\CC{\mathbb{C}}
\def\FF{\mathbb{F}}
\def\PP{\mathbb{P}}
\def\QQ{\mathbb{Q}}
\def\ZZ{\mathbb{Z}}

\def\hbar{\overline{h}}

\newcounter{fixmectr}

\DeclareMathOperator{\an}{an}
\DeclareMathOperator{\dR}{dR}
\DeclareMathOperator{\Frac}{Frac}
\DeclareMathOperator{\MW}{MW}
\DeclareMathOperator{\rig}{rig}
\DeclareMathOperator{\Spec}{Spec}
\DeclareMathOperator{\Trace}{Trace}

\begin{document}

\pagestyle{headings}

\mainmatter

\title{Computing Zeta Functions via $p$-adic Cohomology}

\titlerunning{Computing Zeta Functions via $p$-adic Cohomology}

\author{Kiran S. Kedlaya\thanks{Thanks to Michael Harrison, Joe
Suzuki, and
Fr\'e Vercauteren for helpful comments, and to David Savitt for carefully
reading an early version of this paper.}}

\authorrunning{Kiran S. Kedlaya}

\institute{Department of Mathematics \\
Massachusetts Institute of Technology \\
77 Massachusetts Avenue \\
Cambridge, MA 02139 \\
\email{kedlaya@math.mit.edu} \\
\texttt{http://math.mit.edu/\~{}kedlaya/}}

\maketitle

\begin{abstract}
We survey some recent applications
of $p$-adic cohomology to machine computation
of zeta functions of algebraic varieties over finite fields
of small characteristic, and suggest some new avenues for further exploration.
\end{abstract}

\section{Introduction}

\subsection{The zeta function problem}

For $X$ an algebraic variety over $\FF_q$ (where we write $q = p^n$
for $p$ prime), the zeta function
\[
Z(X,t) = \exp \left( \sum_{i=1}^\infty \frac{t^i}{i} \#X(\FF_{q^i}) \right)
\]
is a rational function of $t$. This fact, the first of the celebrated
Weil Conjectures, follows from Dwork's proof
using $p$-adic analysis \cite{dwork}, 
or from the properties of \'etale ($\ell$-adic) cohomology
(see \cite{freitag-kiehl} for an introduction).

In recent years, the algorithmic problem of determining $Z(X,t)$
from defining equations of $X$ has come into prominence, primarily due to
its relevance in cryptography. Namely,
to perform cryptographic functions using the Jacobian group
of a curve over $\FF_q$, one must first compute the order of said group,
and this is easily retrieved from the zeta function of the curve (as
$Q(1)$, where $Q(t)$ is as defined below).
However, the problem is also connected with other applications of algebraic
curves (e.g., coding theory) and with other computational problems in number
theory (e.g., determining Fourier coefficients of modular forms).

Even if one restricts $X$ to being a curve of genus $g$, in which
case
\[
Z(X,t) = \frac{Q(t)}{(1-t)(1-qt)}
\]
with $Q(t)$ a polynomial over $\ZZ$ of degree $2g$,
there is no algorithm known\footnote{That is, unless one resorts
to quantum computation:
one can imitate Shor's quantum factoring algorithm to
compute 
the order of the Jacobian over $\FF_{q^n}$ for $n$ up to about $2g$, 
and then recover $Z(X,t)$. See \cite{kedlaya-quant}.}
 for computing 
$Z(X,t)$ which is polynomial in the full input size, i.e., in $g$,
$n$, and $\log(p)$. However, if one allows polynomial dependence in
$p$ rather than its logarithm, then one can obtain a polynomial time
algorithm using Dwork's techniques, as shown by Lauder and Wan
\cite{lauder-wan}.
The purpose of this paper is to illustrate how these ideas
can be converted into more practical algorithms in many cases.

This paper has a different purpose in mind than most prior and current
work on computing zeta functions, which has been oriented towards
low-genus curves over large fields (e.g.,
elliptic curves of ``cryptographic size''). This problem is well under
control now; however, we are much less adept at handling curves of high genus
or higher dimensional varieties over small fields. It is in this arena that
$p$-adic methods shoud prove especially valuable; our hope is for this
paper, which mostly surves known algorithmic results on curves, to serve
as a springboard for higher-genus and higher-dimensional investigations.

\subsection{The approach via $p$-adic cohomology}

Historically, although Dwork's proof predated the advent of $\ell$-adic
cohomology, it was soon overtaken as a theoretical tool\footnote{The
gap has been narrowed recently by the work of Berthelot and
others; for instance, in \cite{kedlaya-weil}, one recovers
the Weil conjectures by imitating Deligne's work using $p$-adic
tools.}
by the approach favored by the Grothendieck school, in which context
the Weil conjectures were ultimately resolved by Deligne \cite{deligne}.
The purpose of this paper is to show that by contrast, from
an algorithmic point of view, ``Dworkian'' $p$-adic methods prove to be
much more useful.

A useful analogy is the relationship between
topological and algebraic de Rham cohomology of
varieties over $\CC$. While the topological cohomology is more convenient
for proving basic structural results, computations are often more convenient
in the de Rham setting, since it is so closely linked to defining
equations. The analogy is more than just suggestive: the $p$-adic constructions
we have in mind are variants of and closely related to
algebraic de Rham cohomology, from which they inherit some
computability.

\subsection{Other computational approaches}

There are several other widely used approaches for computing
zeta functions; for completeness, we briefly review these and compare
them with the cohomological point of view.

The method of Schoof \cite{schoof} (studied later by Pila \cite{pila}
and Adleman-Huang \cite{adleman-huang})
is to compute the zeta
function modulo $\ell$ for various small primes $\ell$,
then apply bounds on the coefficients of the zeta function plus the
Chinese remainder theorem.
This loosely corresponds to computing in
$\ell$-adic and not $p$-adic cohomology. This has the benefit of working
well even in large characteristic; on the downside, one can only treat curves,
where $\ell$-adic cohomology can be reinterpreted
in terms of Jacobian varieties, and moreover, one must work with the Jacobians
rather concretely (to extract division polynomials), which is algorithmically
unwieldy. In practice,
Schoof's method has only been deployed in genus 1 (by Schoof's
original work, using improvements
by Atkin, Elkies, Couveignes-Morain, etc.) and genus 2
(by work of Gaudry and Harley \cite{gaudry-harley}, with improvements
by Gaudry and Schost \cite{gaudry-schost}).

A more $p$-adic approach was given by Satoh \cite{satoh}, based on
iteratively computing the Serre-Tate canonical lift \cite{serre-tate}
of an ordinary abelian variety, where one can read off the zeta function
from the action of Frobenius on the tangent space at the origin.
A related idea, due to Mestre, is to compute ``$p$-adic periods''
using a variant of the classical AGM iteration for computing elliptic
integrals. This method has been used to set records for zeta function
computations in characteristic 2 (e.g., \cite{lercier-lubicz}).
The method extends in principle to higher characteristic \cite{kohel}
and genus (see \cite{ritzenthaler1}, \cite{ritzenthaler2} for the genus 3
nonhyperelliptic case),
but it seems difficult to avoid exponential dependence
on genus and practical hangups in handling not-so-small characteristics.

We summarize the comparison between these approaches in the following table.
(The informal comparison in the $n$ column is based on the case
of elliptic curves of a fixed small characteristic.)

\begin{table}[ht] \label{table:comp}
\begin{center}
\begin{tabular}{|c|c|c|c|c|}
\hline
& & \multicolumn{3}{c|}{Dependence on:} \\
Algorithm class & Applicability & $p$ & $n$ & $g$\\
\hline \hline
Schoof & curves & polylog & big polynomial & at least exponential  \\
\hline
Canonical lift/AGM & curves & polynomial & small polynomial &
at least exponential \\
\hline
$p$-adic cohomology & general & 
nearly linear & medium polynomial & polynomial \\
\hline
\end{tabular}
\end{center}
\caption{Comparison of strategies for computing zeta functions}
\end{table}

\section{Some $p$-adic cohomology}

In this section, we briefly describe some constructions of $p$-adic
cohomology, amplifying the earlier remark that it strongly resembles
algebraic de Rham cohomology.

\subsection{Algebraic de Rham cohomology}
We start by recalling how algebraic de Rham cohomology is constructed.
First suppose $X = \Spec A$ is a smooth affine variety\footnote{By ``variety 
over $K$'' we always mean a separated, finite type $K$-scheme.}
over a field $K$
of characteristic zero. Let $\Omega^1_{A/K}$ be the module of
K\"ahler differentials, and put $\Omega^i_{A/K} = \wedge^i \Omega^1_{A/K}$;
these are finitely generated locally free $A$-modules since $X$ is smooth.
By a theorem of Grothendieck \cite{grothendieck}, the cohomology of the complex
$\Omega^i_{A/K}$ is finite dimensional.

If $X$ is smooth but not necessarily affine, one has similar results
on the sheaf level. That is, the hypercohomology of the complex
formed by the sheaves of differentials is finite dimensional. In fact,
Grothendieck proves his theorem first when $X$ is smooth and proper,
where the result follows by a comparison theorem to topological
cohomology (via Serre's GAGA theorem), then uses resolution of singularities
to deduce the general case.

For general $X$, one can no longer use the modules of differentials,
as they fail to be coherent. Instead, following Hartshorne
\cite{hartshorne}, one (locally) embeds $X$ into a smooth scheme $Y$,
and computes de Rham cohomology on the formal completion of $Y$ along $X$.

As one might expect from the above discussion, it is easiest to compute
algebraic de Rham cohomology on a variety $X$ if one is given a good
compactification $\overline{X}$, i.e., a smooth proper variety
such that $\overline{X} \setminus X$ is a normal crossings divisor.
Even absent that, one can still make some headway by computing with
$\mathcal{D}$-modules (where $\mathcal{D}$ is a suitable ring of differential
operators), as shown by Oaku, Takayama, Walther, et al.
(see for instance \cite{walther}).

\subsection{Monsky-Washnitzer cohomology}

We cannot sensibly work with de Rham cohomology directly in characteristic
$p$, because any derivation will kill $p$-th powers and so the cohomology
will not typically be finite dimensional. 
Monsky and Washnitzer \cite{monsky-washnitzer}, \cite{monsky2},
\cite{monsky3} (see also
\cite{vanderput}) introduced a $p$-adic
cohomology which imitates algebraic de Rham cohomology by lifting
the varieties in question to characteristic zero in a careful way.

Let $X = \Spec A$ 
be a smooth affine variety over a finite field $\FF_q$
with $q = p^n$, and let $W$ be the ring of Witt vectors over $\FF_q$,
i.e., the unramified extension of $\ZZ_p$ with residue field $\FF_q$.
By a theorem of Elkik \cite{elkik}, we can find a smooth affine scheme
$\tilde{X}$ over $W$ such that $\tilde{X} \times_W \FF_q \cong X$.
While $\tilde{X}$ is not determined by $X$, we can ``complete along the
special fibre'' to get something more closely bound to $X$.

Write $\tilde{X} = \Spec \tilde{A}$ 
and let $A^\dagger$ be the \emph{weak completion}
of $\tilde{A}$, which is the smallest subring containing
$\tilde{A}$ of the $p$-adic completion of 
$\tilde{A}$ which is $p$-adically saturated (i.e., if $px \in A^\dagger$,
then $x \in A^\dagger$) and closed under the formation of series of the form
\[
\sum_{i_1,\dots,i_m\geq 0} c_{i_1,\dots, i_m} x_1^{i_1}\cdots x_m^{i_m}
\]
with $c_{i_1,\dots,i_m} \in W$ and $x_1, \dots, x_m \in p A^\dagger$.
We call $A^\dagger$ 
the \emph{(integral) dagger algebra} associated to $X$; it is determined by
$X$, but only up to \emph{noncanonical} isomorphism.

In practice, one can describe the weak completion a bit more concretely,
as in the following example.
\begin{lemma}
The weak completion of $W[t_1, \dots, t_n]$ is the ring
$W \langle t_1, \dots, t_n \rangle^\dagger$
of power series over $W$ which converge for
$t_1, \dots, t_n$ within the disc (in the integral closure of $W$)
around $0$ of some radius greater than $1$. 
\end{lemma}
In general, $A^\dagger$
is always a quotient of $W \langle t_1, \dots, t_n \rangle^\dagger$ for some
$n$.

We quickly sketch a proof of this lemma. On one hand,
$W \langle t_1, \dots, t_n \rangle^\dagger$ (which is clearly
$p$-adically saturated) is weakly complete:
if $x_1, \dots, x_m \in pW \langle t_1, \dots, t_n \rangle^\dagger$, then
for $t_1, \dots, t_n$ in some disc of radius strictly greater than 1,
the series defining $x_1, \dots, x_m$ converge to limits of norm 
less than 1,
and so $\sum c_{i_1,\dots, i_m} x_1^{i_1}\cdots x_m^{i_m}$ converges
on the same disc. 
On the other hand, any element of $W \langle t_1, \dots, t_n \rangle^\dagger$
has the form
\[
\sum_{j_1,\dots,j_n \geq 0} e_{j_1,\dots, j_n} t_1^{j_1}\cdots t_n^{j_n}
\]
where $v_p(e_{j_1,\dots,j_n}) + a(j_1+\cdots+j_n) > -b$ for some $a,b$
with $a>0$ (but no uniform choice of $a,b$ is possible). 
We may as well assume that $1/a$ is an integer, and that $b>0$ (since
the weak completion is saturated). Then it is possible to write this series
as
\[
\sum_{i_1,\dots,i_m\geq 0} c_{i_1,\dots, i_m} x_1^{i_1}\cdots x_m^{i_m}
\]
where the $x$'s run over $p^j t_k$ for $j = 1,\dots, 1/a$ and
$k = 1,\dots, n$; hence it lies in the weak completion.

The module of continuous differentials over $A^\dagger$ can be constructed
as follows: given a surjection
$W \langle t_1, \dots, t_n \rangle^\dagger \to A^\dagger$,
$\Omega^1_{A^\dagger}$ is the $A^\dagger$-module 
generated by $dt_1,\dots, dt_n$
modulo enough relations to obtain a well-defined derivation
$d: A^\dagger \to \Omega^1_{A^\dagger}$ satisfying the rule
\[
d \left( \sum_{j_1,\dots,j_n \geq 0} e_{j_1,\dots, j_n}
t_1^{j_1}\cdots t_n^{j_n} \right)
= \sum_{i=1}^n \sum_{j_1,\dots,j_n\geq 0} j_i e_{j_1,\dots, j_n}
(t_1^{j_1}\cdots t_i^{j_i-1}\cdots t_n^{j_n}) dt_i.
\]
Then the Monsky-Washnitzer cohomology (or MW-cohomology) $H^i_{\MW}(X)$ 
of $X$ is the cohomology
of the ``de Rham complex''
\[
\cdots \stackrel{d}{\to} \Omega^i_{A^\dagger} \otimes_W W[\frac 1p]
\stackrel{d}{\to} \cdots,
\]
where $\Omega^i_{A^\dagger} = \wedge^i_{A^\dagger} 
\Omega^i_{A^\dagger}$. Implicit in this
definition is the highly nontrivial fact that this cohomology is 
independent of all of the choices made. Moreover, if $X \to Y$ is a morphism
of $\FF_q$-varieties, and $A^\dagger$ and $B^\dagger$ are corresponding
dagger algebras, then the morphism lifts to a ring map $B^\dagger \to
A^\dagger$, and the induced maps $H^i_{\MW}(Y) \to H^i_{\MW}(X)$
do not depend on the choice of the ring map. The way this works
(see \cite{monsky-washnitzer} for the calculation) is that there are
canonical homotopies (in the homological algebra sense) between any two
such maps, on the level of the de Rham complexes.

MW-cohomology is always finite dimensional over $W[\frac 1p]$; 
this follows from the
analogous statement in rigid cohomology (see \cite{berthelot2}).
Moreover, it admits an analogue of the Lefschetz trace formula for
Frobenius: if $X$ is purely of dimension $d$,
and $F: A^\dagger \to A^\dagger$ is a ring map lifting the
$q$-power Frobenius map, then for all $m > 0$,
\[
\#X(\FF_{q^m}) = \sum_{i=0}^{d} (-1)^i \Trace(q^{dm} F^{-m}, H^i_{\MW}(X)).
\]
This makes it possible in principle, and ultimately in
practice, to compute zeta functions by computing
the action of Frobenius on MW-cohomology.

\subsection{Rigid cohomology}

As in the algebraic de Rham setting, 
it is best to view Monsky-Washnitzer cohomology
in the context of a theory not limited to affine varieties. This context
is provided by Berthelot's rigid cohomology; since we won't compute
directly on this theory, we only describe it briefly. See
\cite{berthelot} or \cite[Chapter~4]{gerkmann}
for a somewhat more detailed introduction.\footnote{We confess that
a presentation at the level of detail we would like does not appear in
print anywhere. Alas, these
proceedings are not the appropriate venue to correct this!}

Suppose $X$ is an $\FF_q$-variety which
is the complement of a divisor in a
smooth proper $Y$ which lifts to
a smooth proper formal $W$-scheme. 
Then this lift gives rise to a rigid analytic
space $Y^{\an}$ via Raynaud's ``generic fibre'' 
construction (its points are the
subschemes of the lift which are integral and finite flat over $W$). This
space comes with a specialization map to $Y$, and the inverse image of $X$
is denoted $]X[$ and called the \emph{tube} of $X$. The rigid
cohomology of $X$ is
the (coherent) cohomology of the direct limit of the de Rham complexes
over all ``strict neighborhoods'' of $]X[$ in $Y^{\an}$. (Within $Y^{\an}$,
$]X[$ is the locus where certain functions take $p$-adic absolute values
less than or equal to 1; to get a strict neighborhood, allow their absolute
values to be less than or equal to $1+\epsilon$ for some $\epsilon>0$.)

For general $X$, we can do the above locally (e.g., on affines) and compute
hypercohomology via the usual spectral sequence; while the construction above
does not sheafify, the complexes involves can be glued ``up to homotopy'',
which is enough to assemble the hypercohomology spectral sequence.

For our purposes, the relevance of rigid cohomology is twofold.
On one hand, it coincides with Monsky-Washnitzer cohomology for $X$
affine. On the other hand, it is related to
algebraic de Rham cohomology via the following theorem.
(This follows, for instance,
from the comparison theorems of \cite{berthelot2} plus the
comparison theorem between crystalline and de Rham cohomology from
\cite{berthelot0}.)
\begin{theorem} \label{thm:compare}
Let $\tilde{Y}$ be a smooth proper $W$-scheme, let $\tilde{Z} \subset 
\tilde{Y}$ be a relative
normal crossings divisor, and set $\tilde{X} = \tilde{Y} \setminus \tilde{Z}$. 
Then there is a canonical isomorphism
\[
H^i_{\dR}(\tilde{X} \times_W (\Frac W)) \to H^i_{\rig}(\tilde{X}
\times_W \FF_q).
\]
\end{theorem}
In particular, if $X$ is affine in this situation,
its Monsky-Washnitzer cohomology is finite dimensional and all of the relations
are explained by relations among algebraic
forms, i.e., relations of finite length. This makes it much easier to
construct ``reduction algorithms'', such as those described in the next
section.

One also has a comparison theorem between rigid cohomology and crystalline
cohomology, a $p$-adic cohomology built in a more ``Grothendieckian''
manner. While crystalline cohomology only behaves well for smooth proper
varieties, it has the virtue of being an \emph{integral} theory. Thus
the comparison to rigid cohomology equips the latter with a canonical
integral structure. By repeating this argument in the context of log-geometry,
one also obtains a canonical integral structure in the setting
of Theorem~\ref{thm:compare}; this is sometimes useful in computations.

\section{Hyperelliptic curves in odd characteristic}

The first\footnote{Although
this seems to be the first overt use of MW-cohomology for 
numerically
computing zeta functions in the literature, it is prefigured by work of 
Kato and Lubkin \cite{kato-lubkin}. Also, similar computations
appear in more theoretical settings, such as Gross's
work on companion forms \cite{gross}.} 
class of varieties where $p$-adic cohomology was demonstrated
to be useful for numerical computations is the class of hyperelliptic curves
in odd characteristic, which we considered in \cite{kedlaya}.
In this section, we summarize the key features of the computation,
which should serve as a prototype for more general considerations.

\subsection{Overview}

An overview of the computation may prove helpful to start with.
The idea is to compute the action of Frobenius on the MW-cohomology of
an affine hyperelliptic curve, and use the Lefschetz trace formula
to recover the zeta function. Of course we cannot compute exactly with
infinite series of $p$-adic numbers, so the computation will be truncated
in both the series and $p$-adic directions, but we arrange to keep
enough precision at the end to uniquely determine the zeta function.

Besides worrying about precision, carrying out this program requires
making algorithmic two features of the Monsky-Washnitzer construction.
\begin{itemize}
\item We must be able to compute a Frobenius lift on a dagger algebra.
\item We must be able to identify differentials forming
a basis of the relevant cohomology
space, and to ``reduce'' an arbitrary differential to a linear combination
of the basis differentials plus an exact differential.
\end{itemize}

\subsection{The dagger algebra and the Frobenius lift}

Suppose that $p \neq 2$,
and let $\overline{X}$ be the hyperelliptic curve of 
genus $g$ given by the affine equation
\[
y^2 = P(x)
\]
with $P(x)$ 
monic of degree $2g+1$ over $\FF_q$ with no repeated roots; in 
particular, $\overline{X}$ has a rational Weierstrass point\footnote{The
case of no rational Weierstrass point is not considered in
\cite{kedlaya}; it has been worked out by Michael Harrison, and has
the same asymptotics.} at infinity.
Let $X$ be the affine curve obtained from $\overline{X}$ by removing all of the
Weierstrass points, i.e., the point at infinity and the zeroes of $y$.

Choose a lift $\tilde{P}(x)$ of $P(x)$
 to a monic polynomial of degree $2g+1$
over $W$. Then the dagger algebra corresponding to $X$ is given by
\[
W \langle x, y, z \rangle^\dagger / (y^2 - \tilde{P}(x), yz - 1),
\]
whose elements can be expressed as $\sum_{i \in \ZZ} A_i(x)y^i$
with $A_i(x) \in W[x]$, $\deg(A_i) \leq 2g$, and $v_p(A_i) + c|i| > d$
for some constants $c,d$ with $c>0$.

The dagger algebra admits a $p$-power Frobenius lift $\sigma$
given by
\begin{align*}
x &\mapsto x^p \\
y &\mapsto y^p \left( 1 + \frac{\tilde{P}(x)^\sigma - 
\tilde{P}(x)^p}{\tilde{P}(x)^p} 
\right)^{1/2},
\end{align*}
which can be computed by a Newton iteration.
Here is where the removal of the Weierstrass points come in handy;
the simple definition of $\sigma$ above clearly requires inverting 
$\tilde{P}(x)$, or
equivalently $y$. It is possible to compute a Frobenius lift on
the dagger algebra of the full affine curve (namely
$W \langle x, y \rangle^\dagger / (y^2 - \tilde{P}(x))$), but this requires
solving for the images of both $x$ and $y$, using a 
cumbersome two-variable Newton iteration.

\subsection{Reduction in cohomology}

The hyperelliptic curve defined by $y^2 = \tilde{P}(x)$,
minus its Weierstrass points, forms a lift $\tilde{X}$
of $X$ of the type described in
Theorem~\ref{thm:compare}, so its algebraic de Rham cohomology coincides
with the MW-cohomology $H^1_{\MW}(X)$.
 That is, the latter is generated
by
\[
\frac{x^i dx}{y} \quad (i=0, \dots, 2g-1), \qquad
\frac{x^i dx}{y^2} \quad (i=0, \dots, 2g)
\]
and it is enough to consider ``algebraic'' relations. Moreover, the cohomology
splits into plus and minus eigenspaces for the hyperelliptic involution
$y \mapsto -y$; the former is essentially the cohomology of $\PP^1$ minus
the images of the Weierstrass points (since one can eliminate $y$ entirely), 
so to compute the zeta function of
$\overline{X}$ we need only worry about the latter. In other words, we
need only consider forms $f(x)dx/y^s$ with $s$ odd.

The key reduction formula is the following: if
$A(x) = \tilde{P}(x) B(x) + \tilde{P}'(x) C(x)$, then
\[
\frac{A(x)\,dx}{y^s} \equiv \left( B(x) + \frac{2C'(x)}{s-2} \right) 
\frac{dx}{y^{s-2}}
\]
as elements of $H^1_{\MW}(X)$. This is an easy consequence of the
evident relation
\[
d \left( \frac{C(x)}{y^{s-2}} \right) \equiv 0
\]
in cohomology.

We use this reduction formula as follows.
Compute
the image under Frobenius of $\frac{x^i dx}{y}$ (truncating large powers of
$y$ or $y^{-1}$, and $p$-adically approximating coefficients).
If the result is
\[
\sum_{j=-M}^N \frac{A_j(x)\,dx}{y^{2j+1}},
\]
use the reduction formula to eliminate the $j=N$ term in cohomology, then
the $j=N-1$ term, and so on, until no terms with $j>0$ remain.
Do likewise with the $j=-M$ term, the $j=-M+1$ term, and so on
(using a similar reduction formula, which
we omit; note that there are relatively few terms on that side anyway).
Repeat for $i=0, \dots, 2g-1$, and construct the ``matrix of the
$p$-power Frobenius'' $\Phi$. Of course the $p$-power Frobenius is not linear,
but the matrix of the $q$-power Frobenius is easily obtained as
$\Phi^{\sigma^{n-1}} \cdots \Phi^\sigma \Phi$, where $\sigma$
here is the Witt vector Frobenius and $q = p^n$.

\subsection{Precision}

We complete the calculation described above
with a $p$-adic approximation of a matrix whose characteristic
polynomial would exactly compute the numerator $Q(t)$ of the zeta function.
However, we can bound the coefficients of that numerator using the Weil
conjectures: if $Q(t) = 1 + a_1 t + \cdots + a_{2g} t^{2g}$,
then for $1 \leq i \leq g$,
$a_{g+i} = q^i a_{g-i}$ and 
\[
|a_i| \leq \binom{2g}{i} q^{i/2}.
\]
In particular, computing $a_i$ modulo a power of $p$ greater than
twice the right side determines it uniquely.

As noted at the end of the previous section, it is critical to know how
much $p$-adic 
precision is lost in various steps of the calculation, in order to
know how much initial precision is needed for the final calculation
to uniquely determine the zeta function. Rather than repeat the whole
analysis here, we simply point out the key estimate
\cite[Lemmas~2 and~3]{kedlaya} and indicate where it comes from.
\begin{lemma} \label{lem:reduction}
For $A_k(x)$ a polynomial over $W$ of degree at most $2g$ and $k \geq 0$
(resp. $k < 0$),
the reduction
of $A_k(x)y^{2k+1}\,dx$ (i.e., the linear combination of $x^i\,dx/y$
over $i=0,\dots,2g-1$ cohomologous to it) becomes integral upon multiplication
by $p^d$ for $d \geq \log_p ((2g+1)(k+1) - 2)$ (resp.
$d \geq \log_p (-2k-1)$).
\end{lemma}
This is seen by considering the polar part of $A_k(x)y^{2k+1}\,dx$
around the point at infinity if $k < 0$, or the other Weierstrass points
if $k > 0$. Multiplying by $p^d$ ensures that the antiderivatives of the
polar parts have integral coefficients, which forces the reductions to do
likewise.

It is also worth pointing out that one can manage precision rather simply
by working in $p$-adic fixed point arithmetic. That is, approximate all numbers
modulo some fixed power of $p$, regardless of their valuation
(in contrast to $p$-adic floating point, where each number is approximated
by a power of $p$ times a mantissa of fixed precision). When a calculation 
produces undetermined high-order digits, fill them in arbitrarily once,
but do not change them later. (That is, if $x$ is computed with some
invented high-order digits, each invocation of $x$ later must use the
\emph{same} invented digits.) The analysis in \cite{kedlaya}, using
the above lemma, shows that most of these invented digits cancel themselves
out later in the calculation, and the precision loss in the reduction process
ends up being negligible compared to the number of digits being retained.

\subsection{Integrality}

In practice, it makes life slightly\footnote{But only slightly: the fact that
there is some basis on which Frobenius acts by an integer matrix means that
the denominators in the product
$\Phi^{\sigma^{n-1}} \cdots \Phi^\sigma \Phi$ can be bounded independently
of $n$.}
easier if one uses a basis in which the matrix
of Frobenius is guaranteed to have $p$-adically \emph{integral} coefficients.
The existence of such a basis is predicted by the comparison with 
crystalline cohomology, but an explicit good basis can be constructed ``by
hand'' by careful use of Lemma~\ref{lem:reduction}. For instance,
the given basis $x^i dx/y$ ($i=0,\dots,2g-1$)
is only good when $p > 2g+1$; on the other hand,
the basis $x^i dx/y^3$ ($i=0,\dots,2g-1$) is good for all $p$ and $g$.

\subsection{Asymptotics}

As for time and memory requirements,
the runtime analysis in \cite{kedlaya} together with
\cite{gaudry-gurel} show that the algorithm requires time
$\tilde{O}(pn^3g^4)$ and space $\tilde{O}(pn^3g^3)$, where again
$g$ is the genus of the curve and $n = \log_p q$. (Here the 
``soft O'' notation ignores logarithmic factors, arising in part
from asymptotically fast integer arithmetic.)

\section{Variations}

In this section, we summarize some of the work on computing MW-cohomology
for other classes of curves. We also mention some experimental results
obtained from implementations of these algorithms.

\subsection{Hyperelliptic curves in characteristic 2}

The method described in the previous section
does not apply in characteristic 2, because the equation
$y^2 = P(x)$ is nonreduced and does not give rise to hyperelliptic curves. 
Instead, one must view the hyperelliptic curve as an Artin-Schreier cover
of $\PP^1$ and handle it accordingly; in particular, we must lift
somewhat carefully. We outline how to do this following
Denef and Vercauteren \cite{denef-vercauteren}, \cite{denef-vercauteren2}.
(Analogous computations based more on Dwork's work have been described by 
Lauder and Wan \cite{lauder-wan2}, \cite{lauder-wan3}, but they seem
less usable in practice.)

Let $\overline{X}$ be a hyperelliptic curve of degree $g$
over $\FF_q$, with $q = 2^n$;
it is defined by some plane equation of the form
\[
y^2 + h(x) y = f(x),
\]
where $f$ is monic of degree $2g+1$ and $\deg(h) \leq g$.
Let $H$ be the monic squarefree polynomial over $\FF_q$ with the same
roots as $h$.
By an appropriate substitution of the form $y \mapsto y + a(x)$,
we can ensure that $f$ vanishes
at each root of $H$.

Let $X$ be the affine curve obtained from $\overline{X}$ by removing
the point at infinity and the zero locus of $H$.
Choose lifts $\tilde{H}, \tilde{h}, \tilde{f}$ 
of $H, h, f$ to polynomials over $W$
of the same degree, such that each root of $\tilde{h}$ is also a root of
$\tilde{H}$, and each root of $\tilde{f}$ 
whose reduction mod $p$ is a root of $H$
is also a root of $\tilde{H}$.
The dagger algebra corresponding to $X$ is now given by
\[
W \langle x, y, z \rangle^\dagger / (y^2 + \tilde{h}(x)y - 
\tilde{f}(x), \tilde{H}(x)z - 1),
\]
and each element can be written uniquely as
\[
\sum_{i \in \ZZ} A_i(x) \tilde{H}_1(x)^i + \sum_{i \in \ZZ} B_i(x) y 
\tilde{H}_1(x)^i
\]
with $\tilde{H}_1(x) = x$ if $\tilde{H}$ is constant and 
$\tilde{H}_1(x) = \tilde{h}(x)$ otherwise,
$\deg(A_i) < \deg(\tilde{H}_1)$ and $\deg(B_i) < \deg(\tilde{H}_1)$ 
for all $i$,
and $v_p(A_i) + c|i| > d$ and $v_p(B_i) + c|i| > d$ for some $c,d$
with $c>0$.
The dagger algebra
admits a Frobenius lift sending $x$ to $x^2$, but this requires
some checking, especially to get an explicit convergence bound; 
see \cite[Lemma~4.4.1]{vercauteren-thesis} for the analysis.

By Theorem~\ref{thm:compare}, the MW-cohomology of $X$
coincides with the cohomology of the hyperelliptic curve $y^2 + \tilde{h}(x) y
- \tilde{f}(x)$ minus the point at infinity and the zero locus of 
$\tilde{H}$. Again,
it decomposes into plus and minus eigenspaces for the hyperelliptic
involution $y \mapsto -y - \tilde{h}(x)$, 
and only the minus eigenspace contributes
to the zeta function of $\overline{X}$. The minus eigenspace is spanned
by $x^i y\,dx$ for $i=0, \dots, 2g-1$, there are again 
simple reduction formulae
for expressing elements of cohomology in terms of this basis,
and one can again bound the precision loss in the reduction; we omit details.

In this case,
the time complexity of the algorithm is $\tilde{O}(n^3 g^5)$ and the space
complexity is $\tilde{O}(n^3 g^4)$. If one restricts to ordinary hyperelliptic
curves (i.e., those where $H$ has degree $g$), the time and space
complexities drop to $\tilde{O}(n^3 g^4)$ and $\tilde{O}(n^3 g^3)$,
respectively, as in the odd characteristic case. It may be possible to 
optimize better for the opposite extreme case, where the curve has
$p$-rank close to zero, but we have not tried to do this.

\subsection{Other curves}

Several variations on the theme developed above have been pursued.
For instance, Gaudry and G\"urel \cite{gaudry-gurel} have
considered superelliptic curves, i.e., those of the form
\[
y^m = P(x)
\]
where $m$ is not divisible by $p$. More generally still,
Denef and Vercauteren consider the class of $C_{a,b}$-curves,
as defined by Miura \cite{miura}. For $a,b$ coprime integers,
a \emph{$C_{a,b}$-curve} is one of the form
\[
y^a + \sum_{i=1}^{a-1} f_i(x) y^i + f_0(x) = 0,
\]
where $\deg f_0 = b$ and $a \deg f_i + bi < ab$ for
$i=1, \dots, a-1$, and the above equation has no singularities in 
the affine plane.

These examples fit into an even broader class of potentially
tractable curves, which we describe following Miura \cite{miura}.
Recall that for a curve $C$ and a point $P$, the
\emph{Weierstrass monoid} is defined to be the set of nonnegative integers
which occur as the pole order at $P$ of some meromorphic function with
no poles away from $P$. Let $a_1 < \dots < a_n$ be a minimal set of generators
of the Weierstrass monoid, and put $d_i = \gcd(a_1, \dots, a_i)$. Then
the monoid is said to be \emph{Gorenstein} (in the terminology
of \cite{nijenhuis-wilf}) if for $i=2, \dots, n$,
\[
\frac{a_i}{d_i} \in \frac{a_1}{d_{i-1}} \ZZ_{\geq 0} 
+ \cdots + \frac{a_{i-1}}{d_{i-1}} \ZZ_{\geq 0}.
\]
If the Weierstrass monoid of $C$ is Gorenstein for some $P$, 
the curve $C$ is said to be \emph{telescopic};
its genus is then equal to
\[
\frac{1}{2} \left( 1 + \sum_{i=1}^n \left( \frac{d_{i-1}}{d_i} - 1 \right)
a_i \right).
\]

The cohomology of telescopic curves is easy to describe, so it seems likely
that one can compute Monsky-Washnitzer cohomology on them. The case $n=2$
is the $C_{a,b}$ case;
for larger $n$, this has been worked out by Suzuki \cite{suzuki}
in what he calls the ``strongly telescopic'' case.
This case is where for each $i$,
the map from $C$ to its image under the projective embedding defined
by $\mathcal{O}(a_i P)$ is a \emph{cyclic} cover (e.g., if $C$ is
superelliptic).

We expect that these can be merged to give an algorithm treating the general
case of telescopic curves.
One practical complication (already appearing in the $C_{a,b}$ case)
is that using a Frobenius lift of the form $x
\mapsto x^p$ necessitates inverting an unpleasantly large polynomial in $y$;
it seems better instead to iteratively compute the action on both $x$ and $y$
of a Frobenius lift without inverting anything.

\subsection{Implementation}

The algorithms described above have proved quite practicable;
here we mention some implementations and report on their performance.
Note that time and space usage figures are only meant to illustrate
feasibility; they are in no way standardized with respect to processor
speed, platform, etc. Also, we believe
all curves and fields described below are 
``random'', without special properties that make them easier to handle.

The first practical test of the original algorithm from \cite{kedlaya} seems
to have been that of Gaudry and G\"urel \cite{gaudry-gurel},
who computed the zeta function of a genus 3 hyperelliptic curve
over $\FF_{3^{37}}$ in 30 hours (apparently not optimized).
They also tested their superelliptic variant, treating a genus 3
curve over $\FF_{2^{53}}$ in 22 hours.

Gaudry and G\"urel \cite{gaudry-gurel2} have also tested the
dependence on $p$ in the hyperelliptic case.
They computed the zeta function of a genus 3 hyperelliptic curve
over $\FF_{251}$ in 42 seconds using 25 MB of memory, and 
over $\FF_{10007}$ in 1.61 hours using 1.4 GB.

In the genus direction, Vercauteren \cite[Sections~4.4--4.5]{vercauteren-thesis}
computed the zeta function of a genus 60 hyperelliptic curve over
$\FF_2$ in 7.64 minutes, and of a 
genus 350 curve over $\FF_2$ in 3.5 days.
We are not aware of any high-genus tests in odd characteristic; in particular,
we do not know whether the lower exponent in the time complexity will
really be reflected in practice.

Vercauteren \cite[Section~5.5]{vercauteren-thesis} has also implemented
the $C_{a,b}$-algorithm in characteristic 2.
He has computed the zeta function of a $C_{3,4}$ curve over $\FF_{2^{288}}$
in 8.4 hours and of a $C_{3,5}$ curve over $\FF_{2^{288}}$
in 12.45 hours.

Finally, we mention an implementation ``coming to a computer near you'':
Michael Harrison has implemented the computation of zeta functions
of hyperelliptic curves in odd characteristic (with or without a rational
Weierstrass point) in a new release of \textsc{Magma}. 
At the time of this writing, we have not seen any performance results.

\section{Beyond hyperelliptic curves}

We conclude by describing some of the rich possibilities for further
productive computations of $p$-adic cohomology, especially in higher
dimensions. A more detailed assessment,
plus some explicit formulae that may prove helpful, appear in the
thesis of Gerkmann \cite{gerkmann} (recently completed under G. Frey).

\subsection{Simple covers}

The main reason the cohomology of hyperelliptic curves in odd 
characteristic is easily computable is that they are ``simple'' (Galois,
cyclic, tamely ramified) covers of a ``simple'' variety (which admits
a simple Frobenius lift). As a first step into higher dimensions, one
can consider similar examples; for instance, a setting
we are currently considering
with de Jong (with an eye toward gathering data on the Tate conjecture
on algebraic cycles)
is the class of double covers of $\PP^2$ of fixed small degree.

One might also consider some simple wildly ramified covers, like
Artin-Schreier covers, which can be treated following Denef-Vercauteren.
(These are also good candidates for Lauder's deformation method; see below.)

\subsection{Toric complete intersections}

Another promising class of varieties to study are smooth complete
intersections in projective space or other toric varieties. These are
promising because their algebraic de Rham cohomology can be computed
by a simple recipe; see \cite[Chapter~5]{gerkmann}.

Moreover, some of these varieties are of current interest thanks to
connections to physics. For instance, Candelas et al. \cite{candelas} 
have studied the zeta functions of
some Calabi-Yau threefolds occurring as toric complete
intersections, motivated by considerations of mirror symmetry.

\subsection{Deformation}

We mention also a promising new technique 
proposed by Lauder.
(A related strategy has been proposed by Nobuo Tsuzuki \cite{tsuzuki}
for computing Kloosterman sums.)
Lauder's strategy is to compute the 
zeta function of a single variety not in isolation, but by placing it
into a family and studying, after Dwork, the variation in Frobenius
along the family as the solution of a certain differential 
equation.\footnote{Lest this strategy seem strangely indirect, note the
resemblance to Deligne's strategy \cite{deligne} for proving the
Riemann hypothesis component of the Weil conjectures!}

A very loose description of the method is as follows. Given an initial 
$X$, say smooth and proper, find a family $f: Y \to B$ over a simple
one-dimensional base (like projective space) which is smooth away from
finitely many points, includes $X$ as one fibre, and has another fibre
which is ``simple''. We also ask for simplicity that the whole situation
lifts to characteristic zero.
For instance, if $X$ is a smooth hypersurface,
$Y$ might be a family which linearly interpolates between the defining
equation of $X$ and that of a diagonal hypersurface. 

One can now compute (on the algebraic lift to characteristic zero)
the Gauss-Manin connection of the family; this will give in particular
a module with
connection over a dagger algebra corresponding to the part of $B$ where
$f$ is smooth. One then shows that there is
a Frobenius structure on this differential equation that computes
the characteristic polynomial of Frobenius on each smooth fibre. That
means the Frobenius structure itself satisfies a differential equation,
which one solves iteratively using an initial condition
provided by the simple fibre. (In the hypersurface example, one can
write down by hand the Frobenius action on the cohomology of a diagonal
hypersurface.)

Lauder describes explicitly
how to carry out the above recipe for Artin-Schreier covers of projective
space \cite{lauder} and smooth projective hypersurfaces \cite{lauder2}. 
The technique has not yet been implemented on a computer, so
it remains to be seen how it performs in practice. It is expected to 
prove most advantageous for 
higher dimensional varieties, as one avoids the need to compute in
multidimensional polynomial rings. In particular, Lauder shows that
in his examples, the dependence of this technique on $d = \dim X$
is exponential in $d$, and not $d^2$. (This is essentially best possible,
as the dimensions of the cohomology spaces in question typically
grow exponentially in $d$.)

\subsection{Additional questions}

We conclude by throwing out some not very well-posed
further questions and suggestions,.

\begin{itemize}
\item Can one can collect data about a class of ``large'' curves
(e.g., hyperelliptic curves of high genus) over a fixed field, and
predict (or even prove) some behavioral properties of the Frobenius
eigenvalues of a typical such curve, in the spirit of Katz-Sarnak?
\item With the help of cohomology computations, can one find nontrivial
instances of cycles on varieties whose existence is predicted by
the Tate conjecture? As noted above, we are looking into this
with Johan de Jong.
\item The cohomology of Deligne-Lusztig varieties furnish 
representations of finite groups of Lie type. Does the $p$-adic cohomology
in particular shed any light on the modular representation theory of these
varieties (i.e., in characteristic equal to that of the underlying field)?
\item There is a close link between $p$-adic Galois representations
and the $p$-adic differential equations arising here;
this is most explicit in the work of Berger \cite{berger}. Can one extend
this analogy to make explicit computations on $p$-adic Galois representations,
e.g., associated to varieties over $\QQ_p$, or modular forms? The work
of Coleman and Iovita \cite{coleman-iovita}
may provide a basis for this.
\end{itemize}

\end{document}